%% file: base-point.tex
\documentclass[12pt]{article}

\usepackage{amsmath}
\usepackage{amssymb}
\usepackage{amsthm}
\usepackage{amsfonts}
\usepackage{amscd}
\usepackage{graphicx}
\usepackage[numbers,sort&compress]{natbib}
\usepackage{longtable}
\usepackage{hyperref}
\usepackage{makeidx}
\usepackage[nottoc]{tocbibind}  
\usepackage[usenames,dvipsnames]{color}
\usepackage{textcomp}
\usepackage{mathrsfs}
\usepackage{array}
\usepackage{marvosym}

\usepackage{enumitem} 

\hypersetup{
    colorlinks=true,
    breaklinks=true,
    filecolor=blue,
    urlcolor=Violet,
    linkcolor=blue,
    citecolor=Green,
    anchorcolor=blue,
    frenchlinks=true,
    pdfborder={0 0 0},
    pdfpagelayout=TwoPageRight,
    pdfdisplaydoctitle=true,
    bookmarksnumbered=true
}

\theoremstyle{plain}

\theoremstyle{definition}
\newtheorem{definition}{Definition}
\newtheorem{remark}{Remark}
\newtheorem{example}{Example}
\newtheorem{algorithm}{Algorithm}

\newcounter{figcnt}

\makeindex
\input{command}

\title{
\vspace{-4.3cm}
Computing basepoints of linear series in the plane
}

\author{Niels Lubbes}
\date{\today}

\begin{document}

\maketitle

\begin{abstract}
We present an algorithm for detecting basepoints of 
linear series of curves in the plane. Moreover, we give an algorithm for constructing
a linear series of curves in the plane for given basepoints.
The underlying method of these algorithms is the classical procedure of blowing up points 
in the plane. We motivate the algorithmic version of this procedure with
several applications.


\end{abstract}

\begingroup
\def\addvspace#1{\vspace{-2mm}}
\tableofcontents
\endgroup

\section{Introduction}

We present two algorithms for analyzing the base locus of a linear series of 
curves in the plane. \ALG{getbp} takes as input a linear series and outputs the 
base locus of this linear series, including also infinitely near basepoints. 
\ALG{setbp} determines which curves in a given linear series form
a linear series with prescribed basepoints.
We resolve basepoints by blowing up points in the plane.
This procedure is classical and well known (see for example \citep[page 28]{har1} or \citep[page 69]{fri1});
however we could not trace back an algorithmic version of this procedure in the literature.
In this article we aim to fill this gap and to advertise this method with some 
applications. See~\citep[\texttt{linear\_series}]{linear_series}
for an open source implementation of 
the algorithms in this article.
The bottleneck of our implementation
is the factorization of univariate polynomials of number fields (\RMK{getbp}). 

Suppose we are given a birational map $\Mdashrow{\McalP}{\MbbP^2}{X}$. 
Rational maps correspond to linear series and
\ALG{getbp} detects the basepoints where $\McalP$ is not defined.
When we know the base locus we can 
compute projective invariants of the rational surface $X$
and with \ALG{setbp} we can compute the linear normalization of $X$.
Moreover, we can
parametrize curves on $X$ that are not 
in the image of $\McalP$. 
This is of interest in computer algebra and geometric modelling~\cite{sen1,sen2}. 
Another application is computing reparametrizations 
and curves on $X$ of given genus and degree~\cite{alc1}.
We will discuss these approaches with examples in \SEC{app}.

\section{Basepoints of linear series}
\label{sec:getbp}

Let $\MbbF$ denote an algebraically closed field of characteristic zero.
Let $F=(f_i)_i$ be a tuple of polynomials \Mst $f_i\in\MbbF[u,v]$.
We think of $F$ as a basis for a linear series of curves in an affine chart 
$\MbbF^2\subset\MbbP^1\times\MbbP^1$ or $\MbbF^2\subset\MbbP^2$.
The curves in this series are of the following form:
\[
C_\alpha(F)=\Mset{ p\in\MbbF^2 }{ \Msum{i}{}\alpha_if_i(p)=0 },
\]
where $\alpha:=(\alpha_i)_i$ with $\alpha_i\in\MbbF$.
We want to resolve the base locus of the map associated to this linear series
using blowups.
In accordance with~\citep[Section I.4]{har1},
the blowup of $\MbbF^2$ at the origin is defined as
\[
U:=\Mset{ (u,v)\times (s:t)}{ ut=vs } \subset\MbbF^2\times\MbbP^1,
\]
and is covered by two charts \Mst $U\cong U_t\cup U_s$, where
\[
U_t:=\{ (\tilde{u},v)\in\MbbF^2 \}\cong\Mset{(u,v)\times (s:t)}{ut=vs,\quad t\neq 0}, 
\]
with $\tilde{u}:=\frac{s}{t}$ \Mst $u=v\tilde{u}$. 
Analogously we define $U_s$ with $s\neq 0$ and $\tilde{v}:=\frac{t}{s}$
\Mst $u\tilde{v}=v$.
The first projection $\Mrow{\pi}{U}{\MbbF^2}$ induces the following maps:
\[
\Mrow{\pi_t}{U_t}{\MbbF^2},\quad(\tilde{u},v)\mapsto (v\tilde{u},v)
\quad\Mand\quad
\Mrow{\pi_s}{U_s}{\MbbF^2},\quad(u,\tilde{v})\mapsto (u,u\tilde{v}).  
\]
We define the translation of the origin to a point $p:=(\hat{u},\hat{v})\in\MbbF^2$:
\[
\Mrow{\tau_p}{\MbbF^2 }{\MbbF^2 },\quad (u,v)\mapsto (u+\hat{u},v+\hat{v}).
\]
Now suppose that 
$p$ is an element in the \Mdef{zero set} $V(F)$
and that~$m>0$
is the multiplicity of $p$ as a point in $C_\alpha(F)$
for general~$\alpha$.
It follows from~\citep[Proposition V.3.6]{har1} that
$C_\alpha(((\tau_p\circ\pi_t)^*f_i)_i)\subset U_t$ 
contains an exceptional component with multiplicity~$m$.
The exceptional component in~$U_t$ 
is where $\pi_t$ is not an isomorphism and thus 
consists of the points with vanishing $v$-coordinate.
It follows that $\gcd(((\tau_p\circ\pi_t)^*f_i)_i)=v^m$.
We define the \Mdef{strict transform} of~$F$ 
\Mwrt basepoint~$p$ and blowup chart~$U_t$ as 
\[
F_{((p,t))}:=((\tau_p\circ\pi_t)^*f_i\div v^m)_i,
\]
where the symbol~$\div$ denotes the polynomial quotient.
Notice that the pullback $(\tau_p\circ\pi_t)^*f_i=f_i(v\tilde{u}+\hat{u},v+\hat{v})$
is a function in $\tilde{u}$ and $v$.

If $V(F_{((p,t))},v)$ or~$V(F_{((p,s))},u)$ is non-empty,
then we can recursively repeat this procedure of computing strict transforms. 
Notice that in the blowup chart~$U_t$, 
we are only interested in basepoints with vanishing $v$-coordinate.
The reason is that the remaining basepoints 
do not lie on the exceptional component and are
already considered in previous charts.

A \Mdef{blowup sequence}~$I$ is defined as either 
an element of $(\MbbF^2\times\{s,t\})^r$ for some $r>0$,
or the empty tuple~$()$.
A strict transform will be denoted $F_I$, where we define $F_{()}$ to be $F$ itself.
For example,~$F_I$ with~$I=((p,s),(q,t))$ denotes
a strict transform along the following sequence of blowups:
we first blowup~$\MbbF^2$ at the basepoint~$p$, 
then we consider the chart~$U_s$,
after that we consider the blowup at the basepoint~$q$,
followed by taking the chart~$U_t$. 
We say that~$q$ is \Mdef{infinitely near} to the basepoint~$p$.

Recall from~\citep[Proposition V.3.2]{har1} that
at each blowup the rank of the Picard group increases. 
The Picard group of a surface that is isomorphic to a blowup of the projective plane, 
is finitely generated. 
It follows that the recursive procedure of computing strict transforms 
for resolving the base locus will always halt.

These considerations lead to the following algorithm.

\begin{algorithm}
\textbf{(get\_basepoints)}
\label{alg:getbp}
\begin{itemize}[topsep=0pt, itemsep=0pt]

\item \textbf{Input:}
A blowup sequence~$I$ and a tuple~$F=(f_i(u,v))_i$ of polynomials~$f_i\in\MbbF[u,v]$
\Mst $\gcd((f_i)_i)\in\MbbF^*$.

\item \textbf{Output:}
A set~$\Gamma$ of elements~$(J,p,m)$ \Mst $J$ is a blowup sequence, 
$p\in V(F_J)$ has multiplicity~$m>0$ and there is 
exactly one strict transform~$F_J$ for each (infinitely near) 
basepoint~$p$.

\item \textbf{Method:}
{\it Let $\frown$ denote concatenation of tuples
and see \SEC{getbp} for the definition of 
$\tau_p$, $\pi_t$ and~$\pi_s$.}

\item[]
If $I=()$, then compute zero set $V(F)$,
else if the last element of the last 2-tuple in $I$ equals $t$,
then compute~$V(F,v)$ 
else compute~$V(F,u)$.
Let $\Upsilon$ denote the computed zero set (see \RMK{getbp}).
Set $\Gamma:=\emptyset$. 

\item[]
For each solution $p\in \Upsilon$ do the following:
\begin{enumerate}
\item[]
$m:=\deg(\gcd(~((\tau_p\circ\pi_t)^*f_i)_i~))$
\item[] 
$\Gamma:=\Gamma~\cup~ \{(I,p,m)\}$
\item[] 
$\Gamma:=\Gamma~\cup~ \textbf{get\_basepoints}(~I^\frown((p,t)),~((\tau_p\circ\pi_t)^*f_i\div v^m)_i~)$
\item[] 
$\Gamma:=\Gamma~\cup~ \textbf{get\_basepoints}(~I^\frown((p,s)),~((\tau_p\circ\pi_s)^*f_i\div u^m)_i~)$
\end{enumerate}

\item[] 
$\textbf{Return}$ $\Gamma$.
\Mend
\end{itemize}
\end{algorithm}

\begin{remark}
\textbf{(implementation)}
\label{rmk:getbp}
\\
The difficulty of implementing \ALG{getbp} is
to compute~$\Upsilon$ over~$\MbbF$. 
This problem can be reduced to the factorization of
a univariate resultant over an algebraic number field~$\MbbK$. 
The roots of nonlinear 
factors of this resultant are adjoined to $\MbbK$ and we factor again, 
until all factors of the resultant are linear. 
The roots of the resultant can be extended to points in $\Upsilon$.
See~\citep[\texttt{linear\_series}]{linear_series} for details about the implementation.
\Mend
\end{remark}

\begin{example}
\textbf{(get\_basepoints)}
\label{exm:getbp}
\\
Let $F=(f_0,f_1):=(u^2 + v^2, v^2 + u)$ and $I:=()$ be input for \ALG{getbp}.
We use the same notation as in the algorithm.
We compute the zero set 
\[
\Upsilon:=V(F)=\{ p_1, p_2, p_3 \},
\]
where $p_1=(0,0)$, $p_2=(1,-\Mi)$ and $p_3=(1,\Mi)$.
We set $\Gamma:=\emptyset$.
Next, we pass into {\it ``for each''}-loop. 
For $p_1\in\Upsilon$ we find that
$((\tau_p\circ\pi_t)^*f_i)_i=( u^2v^2+v^2, v^2+uv  )$
and $\gcd( v^2u^2+v^2, v^2+uv  )=v$ so that $m=1$.
We set
\[
\Gamma:=\Gamma\cup \{\langle I,p_1,1 \rangle\}=\{ \langle (),p_1,1 \rangle \}.
\]
We find that $F_{((p_1,t))}=(f_i\circ\tau_{p_1}\circ\pi_t\div v)_i=(u^2v+v,v+u)$.
Now we call \ALG{getbp} recursively with $F_{((p_1,t))}$ and 
$
\tilde{I}
=()^\frown((p_1,t))
=((p_1,t))
$
and append the output to $\Gamma$: 
\[
\Gamma:=\Gamma\cup
\{ \langle \tilde{I}, p_4, 1\rangle \}=
\{ \langle (),p_1,1 \rangle,~ \langle ((p_1,t)), p_4, 1\rangle \}
,
\]
where $p_4=(0,0)$.
Next we consider the chart~$U_s$ and blowup center~$p_1\in\Upsilon$.
Now $(f_i\circ\tau_{p_1}\circ\pi_s)_i=( u^2 + u^2v^2, u^2v^2 + u  )$
so that~$F_{((p_1,s))}=( u+uv^2, uv+1 )$.
If we call \ALG{getbp} recursively with~$F_{((p_1,s))}$ and~$I:=((p_1,s))$,
then we find that $V( u+uv^2, uv+1, u)=\emptyset$ and thus 
$\Gamma:=\Gamma\cup\{\emptyset\}$.

The algorithm continues with~$p_2,p_3\in\Upsilon$
and the outputs of recursive calls of \ALG{getbp} are again the empty set.
Thus the final output is 
\[
\Gamma=
\{ 
\left\langle ()       , p_1 , 1 \right\rangle,~ 
\left\langle ((p_1,t)), p_4 , 1 \right\rangle,~
\left\langle ()       , p_2 , 1 \right\rangle,~
\left\langle ()       , p_3 , 1 \right\rangle
\}.
\]
Note that this set $\Gamma$ can be expressed as a tree data structure:

$\qquad\Marrow{\emptyset}{}{}$
$\left\langle ()              , p_1, 1 \right\rangle$   
$\Marrow{}{}{}$
$\left\langle ((p_1,t))       , p_4, 1 \right\rangle$

$\qquad\Marrow{\emptyset}{}{}$
$\left\langle ()              , p_2, 1 \right\rangle$

$\qquad\Marrow{\emptyset}{}{}$
$\left\langle ()              , p_3, 1 \right\rangle$

In this case, we say that $F$ has two simple basepoints
and a basepoint of multiplicity one together with an infinitely near basepoint~\citep[page~69]{fri1}.
\Mend
\end{example}

\section{Linear series from basepoints}
\label{sec:setbp}

We present an algorithm which takes as input basepoints $\Gamma$
with some linear series $G$
and outputs the linear series defined by the curves in $G$
that pass through~$\Gamma$.
Typically, $G$ contains all planar curves of some fixed degree.

We have seen in \EXM{getbp} that output~$\Gamma$ of \ALG{getbp}
has the structure of a connected tree \Mst subsequent vertices have nonzero multiplicity. 
We call such~$\Gamma$ a \Mdef{basepoint tree}. 
Before we can state the algorithm for constructing
a linear series with given basepoints,
we introduce notation to access leafs of $\Gamma$.

\begin{definition}
\textbf{(notation for basepoint tree $\Gamma$)} 
\label{def:bptree}
\\
We denote by~$\Gamma_{(I,p,t)}$
the basepoints in the basepoint tree~$\Gamma$ that are in the chart~$U_t$ and that are 
infinitely near to a basepoint~$p\in \MbbF^2$, \Mst $p$ is infinitely near to basepoints specified 
by some blowup sequence~$I$.
Similarly, we denote $\Gamma_{(I,p,s)}$ for the chart~$U_s$.
We denote by $\Gamma_\emptyset$ the leaves at the same level, that are direct neighbors 
of the tree root~$\emptyset$.
\Mend
\end{definition}

If $\Gamma$ is as in \EXM{getbp}, $I=()$ and $J=((p_1,t))$, then 
\[
\Gamma_{(I,p_1,t)}=
\{ 
\langle ((p_1,t)), p_4, 1 \rangle
\},\qquad
\Gamma_{(I,p_1,s)}
=
\Gamma_{(J,p_4,t)}
=
\emptyset,
\]
and $\Gamma_\emptyset
=
\left\{~ 
\langle (), p_1 , 1 \rangle,~
\langle (), p_2 , 1 \rangle,~
\langle (), p_3 , 1 \rangle
~\right\}
$.

\begin{algorithm}
\textbf{(set\_basepoints)}
\label{alg:setbp}
\begin{itemize}[topsep=0pt, itemsep=0pt]
\item \textbf{Input:}
A basepoint tree~$\Gamma$ with nodes~$(I,p,m)$ \Mst $I$ is a blowup sequence,
$p\in\MbbF^2$ and $m>0$.
A tuple~$G=(g_i)_i$ of polynomials~$g_i\in\MbbF[u,v]$.

\item \textbf{Output:}
A matrix~$M$ with entries in~$\MbbF$. 
Let $(k_{ij})_{ij}$ be a column basis of~$\ker M$.
Define
$F=(\sum_ik_{ij}g_i)_j$ 
to be a tuple of polynomials in~$\MbbF[u,v]$.
For all~$(I,p,m)\in \Gamma$, 
the basepoint~$p\in V(F_I)$ has multiplicity~$m>0$
where $F_I$ is a strict transform of~$F$.
Moreover, the base locus of~$G$ is a subscheme
of the base locus of~$F$.
If no $F$ with the above properties exists, then either the empty-list is returned
or the kernel of $M$ is trivial.

\item \textbf{Method:}
{\it See \SEC{getbp} and \DEF{bptree} for the notation of $\pi_t$, $\pi_s$, $\tau_p$, $\Gamma_{(I,p,t)}$ and $\Gamma_{(I,p,s)}$. 
A matrix is represented as a tuple of tuples and $^\frown$ denotes the
appending of rows to a matrix.}

\item[]
$M:=()$
\\[2mm]
For all $(I,p,m)\in\Gamma_\emptyset$ do the following:
\begin{itemize}
\item[] 
$M_0
:=
\left(~ 
\left(
\frac{\partial}{\partial u^a\partial v^b}g_i\circ\tau_p|_{(0,0)}
\right)_i 
~|~
a,b\in\MbbZ_{\geq0},~a+b<m 
~\right)$
\item[] $M_1:=\textbf{set\_basepoints}(~\Gamma_{(I,p,t)},~ ( (\tau_p\circ\pi_t)^*g_i\div v^m)_i ~)$
\item[] $M_2:=\textbf{set\_basepoints}(~\Gamma_{(I,p,s)},~ ( (\tau_p\circ\pi_s)^*g_i\div u^m)_i ~)$
\item[] $M:=M^\frown M_0^\frown M_1^\frown M_2$
\end{itemize}
$\textbf{Return}$ $M$.
\Mend
\end{itemize}

\end{algorithm}

The proof of correctness for \ALG{setbp} 
follows from the observation that this algorithm
recursively traverses through a basepoint tree as \ALG{getbp} would do. 
At each recursion step we add the necessary linear conditions
on~$G$ by adding rows to the matrix~$M$.
Notice that the multiplicity of a basepoint at the origin 
can be expressed by the order of vanishing of partial derivatives.

\begin{example}
\textbf{(set\_basepoints)}
\label{exm:setbp}
\\
Suppose that $\Gamma$ is the output of \ALG{getbp} in \EXM{getbp}:
\[
\Gamma=
\left\{~ 
\left\langle ()       , p_1, 1 \right\rangle,~ 
\left\langle ((p_1,t)), p_4, 1 \right\rangle,~
\left\langle ()       , p_2, 1 \right\rangle,~
\left\langle ()       , p_3, 1 \right\rangle
~\right\}.
\]
We use the same notation as in \ALG{setbp}.
Let $G$ be defined by the monomial basis for quadratic polynomials in $\MbbF[u,v]$:
\[
G=( g_1, g_2, g_3, g_4, g_5, g_6 )=( u^2, uv, u, v^2, v, 1 ). 
\]
We input~$\Gamma$ in \ALG{setbp}. 
We set~$M:=()$ after which we enter the {\it ``for all''}-loop.

We consider the first element $\langle (), p_1, 1 \rangle \in \Gamma_\emptyset$.
The algorithm computes~$M_0$ by evaluating partial derivatives of polynomials in~$G$ 
at~$p_1$. 
Note that~$\tau_{p_1}$ is the identity function.
Since the multiplicity of~$p_1$ is required to be~$1$,
the indices~$(a,b)$ can only attain value~$(0,0)$
so that $\frac{\partial}{\partial u^a\partial v^b}$ is nothing 
but the identity operator.
Thus the matrix~$M_0$ consists of a single row:
\[
M_0=
((
u^2|_{(0,0)}, 
uv|_{(0,0)}, 
u|_{(0,0)}, 
v^2|_{(0,0)}, 
v|_{(0,0)}, 
1|_{(0,0)}
))
=
((0,0,0,0,0,1)).
\]
We recursively call \ALG{setbp} with 
\[
\Gamma_{((),p_1,t)}=
\left\{ 
\left\langle ((p_1,t)), p_4, 1\right\rangle
\right\}
\text{ and }
( 
(\tau_{p_1}\circ\pi_t)^*g_i\div v
)_i
=
(u^2v, uv, u, v,1,0)
\]
and obtain~$M_1=((0,0,0,0,1,0))$ as output.
Next we call \ALG{setbp} with 
\[
\Gamma_{((),p_1,s)}=\{\}
\quad\text{ and }\quad
( 
(\tau_{p_1}\circ\pi_s)^*g_i\div u
)_i
=
(u,uv,1,uv^2,v,0).
 \]
with output~$M_2=()$ so that 
\[
M:=M^\frown M_0^\frown M_1^\frown M_2=((0,0,0,0,0,1),(0,0,0,0,1,0)).
\]
For the second element $\Mmod{(),p_2,1}\in\Gamma_\emptyset$ 
we notice that~$\frac{\partial}{\partial u^0\partial v^0}$ is again the identity
operator and that $\tau_{p_2}$ a translation so that
\[
(\tau_{p_2}^*g_i)_i
=
(
u^2 + 2u + 1, 
uv -\Mi u + v -\Mi, 
u + 1, 
v^2 -2\Mi v - 1, 
v -\Mi, 
1
)
\] 
and thus
$M_0=((1,-\Mi,1,-1,-\Mi,1))$.
We are at a leaf in $\Gamma$ \Mst $M_1=M_2=()$ and thus 
\[
M:=(
(0,0,0,0,0,1),
(0,0,0,0,1,0),
(1,-\Mi,1,-1,-\Mi,1)
).
\]
For the third element $\Mmod{(),p_3,1}\in\Gamma_\emptyset$ 
we find, similarly as before, that
\[
(\tau_{p_3}^*g_i)_i
=
(
u^2 + 2u + 1, 
uv + \Mi u + v + \Mi, 
u + 1, 
v^2 + 2\Mi v - 1, 
v + \Mi, 
1
)
\]
so that
$M_0=((1,\Mi,1,-1,\Mi,1))$ and thus
\[
M:=
(
(0,0,0,0,0,1),
(0,0,0,0,1,0),
(1,-\Mi,1,-1,-\Mi,1),
(1,\Mi,1,-1,\Mi,1)
)
.
\]
We traversed each element of~$\Gamma_\emptyset$  
so the latter~$M$ is the output of \ALG{setbp}.
We are left to compute
\[
\ker M 
= 
\ker
\begin{bmatrix}
0 &    0 & 0 &  0 &    0 & 1 \\ 
0 &    0 & 0 &  0 &    1 & 0 \\
1 & -\Mi & 1 & -1 & -\Mi & 1 \\
1 &  \Mi & 1 & -1 &  \Mi & 1 \\
\end{bmatrix}
\rightsquigarrow
K^\top:=
\begin{bmatrix}
1 & 0 & 0 & 1 & 0 & 0 \\ 
0 & 0 & 1 & 1 & 0 & 0 \\
\end{bmatrix}.
\]
We multiply~$K^\top$ with the column vector defined by~$G$ so that 
$F=(u^2+v^2,u+v^2)$ is the corresponding linear series. 
Indeed this 
was the input of \ALG{getbp} in \EXM{getbp}.
\Mend
\end{example}

\section{Applications}
\label{sec:app}

We propose some examples of applications for \ALG{getbp} and \ALG{setbp}.

\subsection{Parametrization of unreachable curve}
\label{sec:unreachable}

We consider the following birational map:
\begin{equation}
\label{eqn:cubic}
\begin{array}{rrcl}
\McalP\colon & \MbbP^2       & \dashrightarrow & X\subset\MbbP^4 \\
             & (x_0:x_1:x_2) & \longmapsto         & (x_0^2:x_0x_1:x_0x_2:x_1^2:x_1x_2),
\end{array}
\end{equation}
where $X=\Mset{y\in\MbbP^4}{
y_2y_3 - y_1y_4=
y_1y_2 - y_0y_4= 
y_1^2 - y_0y_3=0
}$.
The linear series for $\McalP$ restricted to the affine chart 
$\Mset{x\in\MbbP^2}{x_2\neq 0}$ is $F=(u^2,uv,u,v^2,v)$. 
We apply \ALG{getbp} and we find that the basepoint tree of $F$ is
$\Gamma=\{\langle (), p,1 \rangle\}$ where $p=(0,0)$.
Thus $\McalP$ is not defined at $(0:0:1)$
and we can resolve this locus of indeterminacy by blowing up $\MbbP^2$
with center $(0:0:1)$. 
The resulting blowup is isomorphic to $X$ and the corresponding exceptional curve $E\subset X$
is not reachable by the parametrization $\McalP$. 
In order to parametrize $E$ we consider the strict transform of $F$ along the blowup,
where we use the notation of \SEC{getbp}:
\[
((\tau_p\circ\pi_t)^*f_i\div v)_i=(vu^2,uv,u,v,1). 
\]
The image of $\Mset{(u,v)\in U_t}{v=0}$ via the map
\[
\Marrow{U_t}{}{X},~ (u,v)\mapsto (u^2v:uv:u:v:1)
\]
is a Zariski open set of $E$. 
To cover $E$ completely we would need to compute 
$\Marrow{U_s}{}{X}$ similarly as before.

\subsection{Linear normalization of rational surface}
\label{sec:norm}

Let the surface $Y\subset \MbbP^3$ be the projection of $X\subset\MbbP^4$ with 
center $(1:1:0:0:0)\notin X$ and with birational parametrization
\begin{equation*}
\begin{array}{rrcl}
\McalQ\colon & \MbbP^2       & \dashrightarrow & Y\subset\MbbP^3 \\
             & (x_0:x_1:x_2) & \longmapsto     & (x_0^2-x_0x_1:x_0x_2:x_1^2:x_1x_2),
\end{array}
\end{equation*}
where $Y=\Mset{y\in\MbbP^3}{ y_1^2y_2 - y_1y_2y_3 - y_0y_3^2=0 }$ is singular.
We apply \ALG{getbp} \Mwrt the charts defined by $x_i\neq 0$. 
We find that $\Gamma=\{\langle (), p,1 \rangle\}$ if $x_2\neq 0$ and 
that the basepoint trees are empty for the remaining charts. 
If we apply \ALG{setbp} with input 
$\Gamma$ and a basis $G=(u^2,uv,u,v^2,v,1)$ for quadratic polynomials in $\MbbF[u,v]$,
then we obtain the complete linear series $F=(u^2,uv,u,v^2,v)$. 
After projectivization we recover $\Mdashrow{\McalP}{\MbbP^2}{X}$ from \EQN{cubic},
where $X$ is the linear normalization of $Y$.
Since $X$ is smooth, it is also the projective normalization of $Y$~\citep[II, Example 7.8.4]{har1}.

\subsection{Neron-Severi lattice of rational surface}
\label{sec:ns}

The \Mdef{real structure} of a variety~$X$ is defined as an antiholomorphic involution $\Mrow{\sigma}{X}{X}$.
Thus the real points of~$X$ are the points that are fixed by~$\sigma$.
Maps between real varieties are compatible with the real structure unless explicitly stated otherwise.
A \Mdef{smooth model} of a singular surface~$X$ is a birational morphism $\Marrow{Y}{}{X}$ 
from a nonsingular surface~$Y$, that does not contract exceptional curves.

The Neron-Severi lattice of a surface is well-known concept \citep[page 461]{jha2}.
For our purposes, we give somewhat more explicit description of the data 
associated to this invariant.
The \Mdef{NS-lattice} $N(X)$ of a rational surface~$X\subset \MbbP^n$ 
consists of the following data:
\begin{enumerate}[topsep=0pt, itemsep=0pt]
\item
A unimodular lattice defined by divisor classes on its smooth model~$Y$ modulo numerical equivalence.

\item
A basis for the lattice.
We will consider two different bases for $N(X)$:
\begin{itemize}[topsep=0pt, itemsep=0pt]
\item \Mdef{type 1}: $\Mmod{\Me_0,\Me_1,\ldots,\Me_r}$ 
where the nonzero intersections are\\ $\Me_0^2=1$ and $\Me_j^2=-1$ for $0<j\leq r$,
\item \Mdef{type 2}: $\Mmod{\Ml_0,\Ml_2, \Mp_1,\ldots,\Mp_r}$ 
where the nonzero intersections are\\ $\Ml_0\cdot \Ml_1=1$ and $\Mp_j^2=-1$ for $0<j\leq r$.
\end{itemize}

\item
A unimodular involution $\Mrow{\sigma_*}{N(X)}{N(X)}$ induced by the real structure of~$X$.

\item
A function $\Mrow{h^0}{N(X)}{\MbbZ_{\geq0}}$ assigning the dimension of global sections
of the line bundle associated to a class.

\item
Two distinguished elements $\Mh,\Mk\in N(X)$ corresponding to
class of a hyperplane sections and the canonical class respectively.
\end{enumerate}
We denote the class of a curve~$C\subset X$ in $N(X)$ by $[C]$.
The class~$[F]\in N(X)$ of a linear series~$F$ is defined by the class of any curve
in the linear series.

For example, suppose that $H$ is the linear series 
associated to the identity map $\Marrow{X}{}{X}$, where $X=\MbbP^2$.
In this case 
$N(X)=\Mmod{\Me_0}$ with
$\sigma_*(\Me_0)=\Me_0$, 
$\Mh=[H]=\Me_0$ and $\Mk=-3\Me_0$.
If~$c=\alpha\Me_0$ for~$\alpha\in\MbbZ_{>0}$, then the dimension~$h^0(c)$ equals 
the number of monomial forms of degree~$\alpha$.

If~$H$ is the linear series 
associated to the identity map $\Marrow{X}{}{X}$
\Mst $X=\MbbP^1\times\MbbP^1$,
then
$N(X)=\Mmod{\Ml_0,\Ml_1}$ with  
$\sigma_*=id$,
$\Mh=[H]=\Ml_0+\Ml_1$ and 
$\Mk=-2(\Ml_0+\Ml_1)$.  
If $C\subset X$ is a curve of bidegree~$(\alpha,\beta)\in\MbbZ^2_{\geq0}$, 
then $[C]=\alpha\Ml_0+\beta\Ml_1$.
We can compute $h^0([C])$ with \ALG{setbp} (see \SEC{curve} for an example).

Recall that if $\Mrow{\pi}{X}{X'}$ is the blowup of~$X'$ in a smooth point,
then $N(X)\cong N(X')\oplus \Mmod{e}$ \Mst $e^2=-1$
and $\Mk=\pi^*\Mk'-e$, where $\Mk$ and $\Mk'$ are the canonical classes of $X$
and $X'$ \Mresp.

Suppose that~$H$ is the linear series of some birational map $\Mdashrow{\varphi}{\MbbP^2}{X}$.
With \ALG{getbp} we resolve the base locus of this map:
\begin{equation*}
\begin{array}{r@{}c@{}l}
        &           Z                  &     \\
        &~\swarrow~~~~\searrow\gamma &     \\
\MbbP^2 & \Mdasharrow{}{\varphi}{}      & X          
\end{array}
\end{equation*}
We find that
$N(Z)=\Mmod{\Me_0,\Me_1,\ldots,\Me_r}$ with signature $(-1,1,\ldots,1)$
and  
$\Mk=-3\Me_0+\Me_1+\ldots+\Me_r$
is the canonical class.
Here $\Me_0$ is the class of the pullback of a line in $\MbbP^2$. 
If $i\neq 0$, then $\Me_i$ is the class of the pullback of the exceptional curve 
resulting from the blowup with center~$p_i$.
The pullback of the class of hyperplane sections of~$X$ on $Z$ equals $[H]_{_Z}=\Mh=\alpha_0\Me_0-\alpha_1\Me_1-\ldots-\alpha_r\Me_r$
where $\alpha_0$ is degree of the curves in $H$ and 
$\alpha_i$ is the multiplicity of the curves in~$H$ at basepoint~$p_i$.
Therefore $\sigma_*(\Me_i)=\Me_j$ \Miff $p_i$ is complex conjugate to $p_j$
for $i,j>0$. 

Notice that $Z$ is not always a smooth model for $X$.
In general,
$\gamma=\rho\circ\mu$ where $\Mrow{\mu}{Z}{Y}$ is either the identity map or contracts exceptional curves 
to smooth points and where $\Mrow{\rho}{Y}{X}$ defines the smooth model of $X$.
The NS-lattice $N(X)$ is in this case recovered as the sublattice of $N(Z)$ that is
orthogonal to the set
$\Mset{c\in N(Z)}{~ h^0(c)=1,~ [H]_{_Z}\cdot c=0,~ c^2=-1 ~}$.
In practice, it is easier to do computations directly in $N(Z)$.
For example, with \ALG{setbp} we can compute 
the dimension of global sections~$h^0(c)$ for any~$c\in N(Z)$
and $h^0([C]_{_Z})=h^0([C])$ for any curve~$C\subset X$, where 
$[C]_{_Z} \in N(Z)$ is the class of the pullback of $C$ to $Z$.
Moreover, the degree of~$C\subset X$ equals $[H]_{_Z}\cdot [C]_{_Z}$.

\subsection{Projective invariants of rational surface}
\label{sec:inv}

We show via examples that we can compute with 
\ALG{getbp} and \ALG{setbp} the following projective invariants of the
rational surface $X$ and linear series $F$ as defined at \EQN{cubic}.
\begin{itemize}[topsep=0pt,itemsep=0pt]
\item \textbf{NS-lattice of~$X$.}
It follows from \SEC{ns} that
$N(X)=\Mmod{\Me_0,\Me_1}$ with 
$\sigma_*=id$,
$\Mh=[F]=2\Me_0-\Me_1$ and 
$\Mk=-3\Me_0+\Me_1$.

\item \textbf{Degree of~$X$.} 
The curves~$C_\alpha(F)$ 
in the linear series are pullbacks of hyperplane sections and thus 
the number of intersections $C_\alpha(F)\cap C_\beta(F)$ for generic 
$\alpha,\beta\in\MbbP^4$ outside the basepoints, equals the degree of~$X$. 
It follows that $\deg X=[F]^2=(2\Me_0-\Me_1)^2=3$.

\item \textbf{Sectional genus~$p_s(X)$.}
We recall that the sectional genus of the surface~$X$ is defined as the geometric
genus of its general hyperplane section:
$p_s(X)=p_a [F]=\frac{1}{2}([F]^2+[F]\cdot \Mk)+1=5$ where~$[F]\cdot \Mk=5$.

\item \textbf{Arithmetic genus~$p_a(X)$.}
It follows from Riemann-Roch theorem and Kodaira vanishing that 
$h^0([F])=\frac{1}{2}([F]^2-[F]\cdot \Mk)+p_a(X)+1$.
Here~$h^0([F])$ corresponds to the number of generators of~$F$
in case~$F$ is a complete linear series.
Recall from \SEC{norm} that we can compute the missing generators 
of incomplete linear series.
With~$X$ as in \EQN{cubic} we find that $h^0([F])=5$ 
and thus $p_a(X)=5$.

\end{itemize}

\subsection{Adjoint surface}

In this section we compute for a given rational surface its adjoint surface.
See \citep{nls-f6} for an application of adjoint surfaces.
Suppose that $H$ is a linear series corresponding to 
the following birational map:
\begin{equation*}
\begin{array}{rrcl}
\McalR\colon & \MbbP^2       & \dashrightarrow & Z\subset\MbbP^{16} \\
             & (x_0:x_1:x_2) & \longmapsto     & 
(
x_1^5:x_1^4x_2:x_1^4x_0:x_1^3x_2^2:x_1^3x_2x_0:x_1^3x_0^2:x_1^2x_2^3: \\&&& 
x_1^2x_2^2x_0:x_1^2x_2x_0^2:x_1^2x_0^3:x_1x_2^4:x_1x_2^3x_0:x_1x_2^2x_0^2:      \\&&& 
x_1x_2x_0^3:x_2^5 - x_2^2x_0^3:x_2^4x_0-x_2^2x_0^3: x_2^3x_0^2 - x_2^2x_0^3
).
\end{array}
\end{equation*}
We apply \ALG{getbp} and find that in the chart~$x_0\neq 0$ that 
\[
\Gamma_Z=\{\langle ((),p_1,2)\rangle, \langle ((),p_2,1)\rangle\},
\]
where $p_1=(0,0)$ and $p_2=(0,1)$.
There are no basepoints outside this affine chart of~$\MbbP^2$. 
We find that 
$N(Z)=\Mmod{\Me_0,\Me_1,\Me_2}$ with $\sigma_{Z*}=id$, 
$\Mk_Z=-3\Me_0+\Me_1+\Me_2$ and $\Mh_Z=[H]=5\Me_0-2\Me_1-\Me_2$,
where~$\Me_i$ is the class of the pullback of the exceptional curve along the blowups with centers~$p_i$.
The parametrization of the adjoint surface has class $\Mh_Z+\Mk_Z=2\Me_0-\Me_1$. 
Thus the adjoint surface of~$Z$ is~$X$ as defined at \EQN{cubic}.

\subsection{Curves on rational surfaces and reparametrizations}
\label{sec:curve}

We start by constructing a sextic del Pezzo surface $S\subset\MbbP^6$ as the blowup
of~$\MbbP^1\times\MbbP^1$ in complex conjugate points~\mbox{$P_1=(1:\Mi;1:-\Mi)$} 
and~\mbox{$P_2=(1:-\Mi;1:\Mi)$}.
The forms of bi-degree~$(2,2)$ on~$\MbbP^1\times\MbbP^1$ are 
\[
(x_0^2y_0^2,~ x_0^2y_0y_1,~ x_0^2y_1^2,~ x_0x_1y_0^2,~ x_0x_1y_0y_1,~ x_0x_1y_1^2,~ x_1^2y_0^2,~ x_1^2y_0y_1,~ x_1^2y_1^2).
\]
We consider the chart of~$\MbbP^1\times\MbbP^1$ where $x_0,y_0\neq 0$ and set 
\[
G_S=(1,v,v^2,u, uv, uv^2, u^2, u^2v, u^2v^2)
\]
as the corresponding dehomogenization of the above forms.
We call \ALG{setbp} with~$G_S$ and
$\Gamma_S=\{\langle p_1,1\rangle,~\langle p_2,1\rangle\}$
where $p_1=(\Mi,-\Mi)$ and $p_2=(-\Mi,\Mi)$.
The output is linear series~$L$, which after bi-homogenization 
defines a birational map:
\begin{equation*}
\begin{array}{rrcl}
\McalS\colon & \MbbP^1\times\MbbP^1 & \dashrightarrow & S\subset\MbbP^6 \\
             & (x_0:x_1;y_0:y_1)    & \longmapsto     & 
(
x_0^2y_0^2-x_1^2y_1^2:
x_0^2y_0y_1+x_1^2y_0y_1:
x_0^2y_1^2+x_1^2y_1^2:\\&&&
x_0x_1y_0^2-x_1^2y_0y_1:
x_0x_1y_0y_1-x_1^2y_1^2:\\&&&
x_1^2y_0y_1+x_0x_1y_1^2:
x_1^2y_0^2+x_1^2y_1^2
).
\end{array}
\end{equation*}
Indeed, $S$ is isomorphic to the blowup of $\MbbP^1\times\MbbP^1$ at~$P_1$ and~$P_2$.
As in \SEC{inv}
we find that $N(S)=\Mmod{\Ml_0,\Ml_1,\Mp_1,\Mp_2}$ with 
$\sigma_*(\Ml_0)=\Ml_0$, 
$\sigma_*(\Ml_1)=\Ml_1$,
$\sigma_*(\Mp_1)=\Mp_2$ 
and
$\Mh_S=-\Mk_S=[L]=2\Ml_0+2\Ml_1-\Mp_1-\Mp_2$.

The images $\McalS(\{\alpha\}\times\MbbP^1)$ 
and $\McalS(\MbbP^1\times\{\alpha\})$ for $\alpha\in\MbbP^1$,
are \Mresp~a conic in the first- and the second- family of conics that cover $S$.
The classes in $N(S)$ of conics in these families are $\Ml_0$ and $\Ml_1$, \Mresp.
The degree of a curve $C\subset S$ is equal to $[L]\cdot [C]$.
There is a third family of conics in $S$ that has class $\Ml_0+\Ml_1-\Mp_1-\Mp_2\in N(S)$.
We construct with \ALG{setbp} a linear series of forms of bi-degree $(1,1)$ with simple basepoints at $P_1$ and $P_2$.
Thus the input is $\Gamma_S$ with $G_S'=(1, v, u, uv)$ and we call the output linear series $L'$ 
where $[L']=\Ml_0+\Ml_1-\Mp_1-\Mp_2$.
After bi-homogenization, $L'$ is equal to $(x_0y_0 - x_1y_1, x_1y_0 + x_0y_1)$.
Let
\[
B_\alpha=\Mset{ (x_0:x_1;y_0:y_1)\in\MbbP^1\times\MbbP^1 }{\alpha_0 (x_0y_0 - x_1y_1) = \alpha_1 (x_1y_0 + x_0y_1) }.
\]
The image~$\McalS(B_\alpha)$ 
defines a conic in the third family, for all~$\alpha\in \MbbP^1$.
Moreover, we verified by computation that $h^0(\Ml_0+\Ml_1-\Mp_1-\Mp_2)=2$.
After rearranging terms, the defining equation of 
$B_\alpha$ is $y_0(x_0\alpha_0-x_1\alpha_1)=y_1(x_0\alpha_1+x_1\alpha_0)$.
It follows that the following birational map parametrizes all curves~$B_\alpha$ 
as $\alpha$ varies in~$\MbbP^1$:
\begin{equation*}
\begin{array}{rrcl}
\McalT\colon & \MbbP^1\times\MbbP^1        & \dashrightarrow & \MbbP^1\times\MbbP^1 \\
             & (x_0:x_1;\alpha_0:\alpha_1) & \longmapsto     & (x_0:x_1:x_0\alpha_1+x_1\alpha_0:x_0\alpha_0-x_1\alpha_1).
\end{array}
\end{equation*}
The composition $\McalS\circ\McalT$ is a reparametrization of~$S$
so that $(\McalS\circ\McalT)(\MbbP^1\times\{\alpha\})$
for $\alpha\in\MbbP^1$,
is now a conic in the third family of conics.

\section{Acknowledgements}

I would like to thank Josef Schicho. 
His many valuable comments and explanations concerning
the content of this paper cannot be overestimated.
I would like to thank Matteo Gallet for his suggestions
which greatly improved the presentation of this paper.
We used \citep[Sage]{sage} for the computations.

\bibliography{geometry}
\paragraph{address of author:}
Johann Radon Institute for Computational and Applied 
Mathematics (RICAM), Austrian Academy of Sciences
\\
\textbf{email:} niels.lubbes@gmail.com

\end{document}

%% file: command.tex
    \newcommand{\Mend}{\hfill \ensuremath{\vartriangleleft}}

    \newcommand{\Msum}[2]{\ensuremath{\overset{#2}{\underset{#1}{\sum}}}}

    \newcommand{\Mmod}[1]{\langle #1\rangle} 

    \newcommand{\Mset}[2]{\ensuremath{\{~ #1 ~|~ #2 ~\}}}

    \newcommand{\Mand}{{\textrm{~and~}}}

    \newcommand{\Miff}{if and only if }

    \newcommand{\Mst}{{such that }}
    
    \newcommand{\Mresp}{respectively}
    \newcommand{\Mwrt}{with respect to }

    \newcommand{\Marrow}[3]{\ensuremath{#1\stackrel{#2}{\longrightarrow}#3}}
    
    \newcommand{\Mdasharrow}[3]{\ensuremath{#1\stackrel{#2}{\dashrightarrow}#3}}
    
    \newcommand{\Mrow}[3]{\ensuremath{#1\colon #2\longrightarrow #3}}
    \newcommand{\Mdashrow}[3]{\ensuremath{#1\colon #2\dashrightarrow #3}}
    
    \newcommand{\Mdef}[1]{\textit{#1}\index{#1}}

    \newcommand{\Me}{e}
    \newcommand{\Mp}{\varepsilon}
    \newcommand{\Mk}{k}
    \newcommand{\Ml}{\ell}
    \newcommand{\Mh}{h}
    \newcommand{\Mi}{\mathfrak{i}}
    
    \newcommand{\McalP}{{\mathcal{P}}}\newcommand{\McalQ}{{\mathcal{Q}}}\newcommand{\McalR}{{\mathcal{R}}}\newcommand{\McalS}{{\mathcal{S}}}\newcommand{\McalT}{{\mathcal{T}}}
    \newcommand{\MbbF}{{\mathbb{F}}}\newcommand{\MbbK}{{\mathbb{K}}}\newcommand{\MbbP}{{\mathbb{P}}}\newcommand{\MbbZ}{{\mathbb{Z}}}

    \newcommand{\DEF}[1]{{Definition~\ref{def:#1}}}
    \newcommand{\RMK}[1]{{Remark~\ref{rmk:#1}}}
    \newcommand{\EQN}[1]{{(\ref{eqn:#1})}}
    \newcommand{\SEC}[1]{{\textsection\ref{sec:#1}}}
    
    \newcommand{\ALG}[1]{{Algorithm~\ref{alg:#1}}}
    
    \newcommand{\EXM}[1]{{Example~\ref{exm:#1}}}